\documentclass[11pt]{amsart}
\usepackage{amsmath,amsthm,amsfonts,amssymb,mathrsfs, mathtools}
\usepackage{enumerate}
\usepackage[colorlinks=true, linkcolor=blue, citecolor=magenta, menucolor=black]{hyperref}
\usepackage[demo]{graphicx}
\usepackage[up]{caption}
\usepackage{subcaption}
\usepackage{pgf,tikz}
\usepackage[toc,page]{appendix}
\usetikzlibrary{arrows}
\usetikzlibrary{patterns}
\usepackage{color}

\numberwithin{equation}{section}
\theoremstyle{plain}

\newtheorem{letterthm}{Theorem}

\newtheorem{lettercor}[letterthm]{Corollary}

\theoremstyle{definition}

\newtheorem*{defn*}{Definition}

\newtheorem*{remark}{Remark}

\newcommand{\N}{\mathbb{N}}
\newcommand{\R}{\mathbb{R}}
\newcommand{\C}{\mathbb{C}}
\newcommand{\Z}{\mathbb{Z}}

\newcommand{\dpr}{^{\prime\prime}}

\newcommand{\rB}{\operatorname{B}}

\newcommand{\rE}{\operatorname{ E}}

\newcommand{\rL}{\operatorname{ L}}

\allowdisplaybreaks

\makeatletter
\@namedef{subjclassname@2020}{%
  \textup{2020} Mathematics Subject Classification}
\makeatother

\begin{document}

\title[Selfless W$^*$-probability spaces]{Selfless W$^*$-probability spaces and \\ Connes'\! bicentralizer problem}

\begin{abstract}
We introduce the notion of selfless W$^*$-probability space and study its connection with Connes'\! bicentralizer problem. In particular, we show that if $M$ is a separable type ${\rm III_1}$ factor with trivial bicentralizer, then $(M, \varphi)$ is selfless for every faithful normal state $\varphi \in M_\ast$.
\end{abstract}

\author{Cyril Houdayer}
\address{\'Ecole normale sup\'erieure \\ D\'epartement de math\'ematiques et applications \\ Universit\'e Paris-Saclay \\ 45 rue d'Ulm \\ 75230 Paris Cedex 05 \\ France}
\email{cyril.houdayer@ens.psl.eu}
\thanks{CH is supported by ERC Advanced Grant NET 101141693}

\author{Amine Marrakchi}
\address{CNRS \\ \'Ecole normale sup\'erieure de Lyon \\ Unit\'e de math\'ematiques pures et appliquées  \\ 46 all\'ee d’Italie \\ 69364 Lyon \\ France}
\email{amine.marrakchi@ens-lyon.fr}


\maketitle

\section{Introduction}

In \cite{Ro23}, Robert introduced a new class of C$^*$-probability spaces, which he called \emph{selfless}, characterized by the existence of a copy of themselves in their ultrapower that is freely independent from the diagonal copy (thus being ``free from themselves"). This property quickly attracted the attention of numerous researchers as it implies many important regularity properties and is satisfied by a large class of examples (see \cite{AGKEP24, HKER25, RTV25, Vi25, Oz25, FKOCP25, Vi26}).

In this short note, we introduce a parallel notion of selfless W$^*$-probability space and we relate this notion to Connes'\! bicentralizer problem.

A W$^*$-{\em probability space} is a pair $(M, \varphi)$ that consists of a von Neumann algebra $M$ endowed with a faithful normal state $\varphi \in M_\ast$. For W$^*$-probability spaces $(M, \varphi)$ and $(N, \psi)$, we say that $(M, \varphi) \subset (N, \psi)$ is an {\em inclusion} of W$^*$-probability spaces if $M \subset N$ and if there exists a faithful normal conditional expectation $\rE : N \to M$ such that $\varphi\circ \rE = \psi$. In that case, $\rE : N \to M$ is the unique faithful normal conditional expectation such that $\varphi \circ \rE = \psi$. 

Following \cite{GH21}, we say that an inclusion of W$^*$-probability spaces $(M, \varphi) \subset (N, \psi)$ is {\em existentially closed} if there exists a nonprincipal ultrafilter $\mathcal U$ on some directed set $I$ such that $(M, \varphi) \subset (N, \psi) \subset (M, \varphi)^{\mathcal U}$, where $(M, \varphi) \subset (M, \varphi)^{\mathcal U}$ is the diagonal inclusion. Note that if  $N$ is separable (i.e.\! $N$ has separable predual), then $\mathcal U$ can be chosen to be a nonprincipal ultrafilter on $\N$.

Adapting \cite[Definition 2.1]{Ro23} to the von Neumann algebraic realm, we say that a W$^*$-probability space $(M, \varphi)$ is {\em selfless} if the first factor inclusion $(M, \varphi) \subset (M, \varphi) \ast (M, \varphi)$ is existentially closed. 

Popa's seminal work \cite{Po95} shows that $(M, \tau)$ is selfless for any separable type ${\rm II_1}$ factor $M$ endowed with its canonical trace $\tau$. Houdayer--Isono \cite{HI14} extended Popa's result by showing that $(M, \varphi)$ is selfless for any separable factor $M$ endowed with a faithful normal state $\varphi \in M_\ast$ for which $(M_\varphi)' \cap M = \C1$.

In this note, we show that a diffuse separable W$^*$-probability space $(M, \varphi)$ is selfless if and only if it has a trivial \emph{bicentralizer}. Recall from \cite{Co80, Ha85}, that the {\em bicentralizer} $\rB(M, \varphi)$ of a W$^*$-probability space $(M, \varphi)$ is the set of all elements $x \in M$ that satisfy the following condition: 
\begin{verse}
\noindent
For every $\varepsilon > 0$, there exists $\delta > 0$ such that for every $u \in \mathscr U(M)$, if $\|u \varphi - \varphi u\| < \delta$, then $\|u x - x u\|_\varphi < \varepsilon$. 
\end{verse}

We then obtain the following characterization.

\begin{letterthm} \label{thm-characterization}
Let $(M, \varphi)$ be a diffuse separable $\mathrm{W}^*$-probability space. The following assertions are equivalent:
\begin{enumerate}[{\rm (i)}]
\item $(M, \varphi)$ is selfless.

\item The first factor inclusion $(M, \varphi) \subset (M, \varphi) \ast (N, \psi)$ is existentially closed for some nontrivial $\mathrm{W}^*$-probability space $(N, \psi)$.

\item The first factor inclusion $(M, \varphi) \subset (M, \varphi) \ast (\rL(\Z), \tau_\Z)$ is existentially closed, where $\tau_\Z$ is the canonical trace on $\rL(\Z)$. 

\item $\rB(M, \varphi) = \C 1$.
\end{enumerate}
\end{letterthm}

We point out that the first three conditions above are analogous to the ones appearing in \cite[Theorem 2.6]{Ro23} for C$^*$-probability spaces.

The bicentralizer $\rB(M, \varphi) \subset M$ is a von Neumann subalgebra such that $\rB(M, \varphi) \subset (M_\varphi)' \cap M$ (see \cite[Proposition 1.3]{Ha85}). Thus Theorem \ref{thm-characterization} strengthens Houdayer--Isono's result \cite[Theorem A]{HI14}. In fact, the proof of Theorem \ref{thm-characterization} relies on \cite[Theorem A]{HI14} combined with a diagonal argument.

Observe that if $\rB(M, \varphi) = \C1$, then $M$ must be a factor. Moreover, according to \cite{Ok21}, one and exactly one of the following assertions hold:
\begin{itemize}
\item $(M, \varphi)$ is a tracial factor of type ${\rm I}_n$ for $n \in \N^*$ or of type ${\rm II_1}$.
\item There exists $\lambda \in (0, 1)$ such that $(M, \varphi)$ is a type ${\rm III_\lambda}$ factor endowed with its $\frac{2\pi}{|\log(\lambda)|}$-periodic faithful normal state. In that case, we have $(M_\varphi)' \cap M = \C 1$.
\item $M$ is a type ${\rm III_1}$ factor. In that case, using \cite[Corollary 1.5]{Ha85}, we further have $\rB(M, \psi) = \C1$ for \emph{every} faithful normal state $\psi \in M_\ast$.
\end{itemize}

In particular, Theorem \ref{thm-characterization} implies that no W$^*$-probability space of type ${\rm II_\infty}$ or type ${\rm III_0}$ can be selfless (this also follows by combining \cite[Theorem 3.5]{GH21} and \cite[Theorem 4.1]{Ue11}).

On the other hand, a famous conjecture of Connes, known as Connes'\! bicentralizer problem, claims that $\rB(M,\varphi)=\C 1$ for \emph{every} type ${\rm III_1}$ factor $M$ and every faithful normal state $\varphi \in M_*$.  This conjecture has been verified for several families of type ${\rm III_1}$ factors such as amenable factors \cite{Ha85}, factors with a Cartan subalgebra, free products \cite{HU15}, semisolid factors \cite{HI15}, $q$-deformed Araki--Woods factors \cite{HI20, Bi24} and tensor products of type ${\rm III_1}$ factors \cite{Ma25}.

\begin{lettercor}\label{cor}
Let $M$ be a separable type ${\rm III_1}$ factor satisfying Connes'\! bicentralizer conjecture. Then for {\em every} faithful normal state $\varphi \in M_\ast$, the $\mathrm{W}^*$-probability space $(M, \varphi)$ is selfless.
\end{lettercor}


\section{Proof of Theorem \ref{thm-characterization}}

\begin{proof}

$(\rm i) \Rightarrow (\rm ii)$ This is obvious by taking $(N, \psi) = (M, \varphi)$.

$(\rm ii) \Rightarrow (\rm iii)$ Let $(N, \psi)$ be a nontrivial W$^*$-probability space such that the first factor inclusion $(M, \varphi) \subset (M, \varphi) \ast (N, \psi)$ is existentially closed. We may assume that $N$ is separable. Choose a nonprincipal ultrafilter $\mathcal U$ on $\N$ such that we have $(M, \varphi) \subset (M, \varphi) \ast (N, \psi) \subset (M, \varphi)^{\mathcal U}$, where $(M, \varphi) \subset (M, \varphi)^{\mathcal U}$ is the diagonal inclusion. There are two cases to consider. 

Firstly, assume that $N_\psi = \C 1$. Then $N$ is a type ${\rm III_1}$ factor (see e.g.\! \cite[Lemma 5.3]{AH12}). Consider the same nonprincipal ultrafilter $\mathcal V = \mathcal U$ on $\N$ and define $\mathcal W = \mathcal V \otimes \mathcal U = \mathcal U \otimes \mathcal U$, which is a nonprincipal ultrafilter on $\N \times \N$. By \cite[Proposition 2.5]{AHHM18}, we naturally have $(M^{\mathcal U}, \varphi^{\mathcal U})^{\mathcal V} = (M^{\mathcal U}, \varphi^{\mathcal U})^{\mathcal U} = (M, \varphi)^{\mathcal W}$. Then we have 
$$(M, \varphi) \subset (M, \varphi) \ast (N, \psi)^{\mathcal U} \subset (M, \varphi)^{\mathcal U} \ast (N, \psi)^{\mathcal U}  \subset (M^{\mathcal U}, \varphi^{\mathcal U})^{\mathcal U} = (M, \varphi)^{\mathcal W},$$ where $(M, \varphi) \subset (M, \varphi)^{\mathcal W}$ is the diagonal inclusion. By \cite[Proposition 4.24]{AH12}, the centralizer $(N^{\mathcal U})_{\psi^{\mathcal U}}$ is a type ${\rm II_1}$ factor and so $(\rL(\Z), \tau_{\Z}) \subset ((N^{\mathcal U})_{\psi^{\mathcal U}}, \psi^{\mathcal U})$. This futher implies that 
$$(M, \varphi) \subset (M, \varphi) \ast (\rL(\Z), \tau_\Z) \subset (M, \varphi) \ast ((N^{\mathcal U})_{\psi^{\mathcal U}}, \psi^{\mathcal U} )  \subset (M, \varphi)^{\mathcal W},$$ where $(M, \varphi) \subset (M, \varphi)^{\mathcal W}$ is the diagonal inclusion. Therefore, the first factor inclusion $(M, \varphi) \subset (M, \varphi) \ast (\rL(\Z), \tau_\Z)$ is existentially closed.

Secondly, assume that $N_\psi \neq \C 1$. Upon replacing $(N, \psi)$ by $(N_\psi, \psi)$, we may assume that $(N, \psi)$ is tracial. Reasoning as is the first case, we obtain
\begin{align*}
(M, \varphi) &\subset (M, \varphi) \ast (N, \psi)^{\ast 2} \\
&= \left( (M, \varphi) \ast (N, \psi) \right ) \ast (N, \psi) \\ 
&\subset (M, \varphi)^{\mathcal U} \ast (N, \psi) \\
&\subset  (M, \varphi)^{\mathcal U} \ast (N, \psi)^{\mathcal U} \\
& \subset (M^{\mathcal U}, \varphi^{\mathcal U})^{\mathcal U} = (M, \varphi)^{\mathcal W}.
\end{align*}
By \cite[Lemma 2.5]{Ro23}, there exists $n \in \N$ large enough so that the iterated free product $(N, \psi)^{\ast n}$ is diffuse and so $(\rL(\Z), \tau_{\Z}) \subset (N, \psi)^{\ast n}$. Upon iterating $n$ times the ultraproduct construction and replacing $\mathcal U$ by the appropriate ultrafilter $\mathcal W = \mathcal U^{\otimes n}$, we obtain 
$$(M, \varphi) \subset (M, \varphi) \ast (\rL(\Z), \tau_\Z) \subset (M, \varphi) \ast (N, \psi)^{\ast n} \subset (M, \varphi)^{\mathcal W},$$ where $(M, \varphi) \subset (M, \varphi)^{\mathcal W}$ is the diagonal inclusion. Therefore, the first factor inclusion $(M, \varphi) \subset (M, \varphi) \ast (\rL(\Z), \tau_\Z)$ is existentially closed.

$(\rm iii) \Rightarrow (\rm iv)$ Choose a nonprincipal ultrafilter $\mathcal U$ on $\N$ such that $(M, \varphi) \subset (M, \varphi) \ast (\rL(\Z), \tau_\Z) \subset (M, \varphi)^{\mathcal U}$, where $(M, \varphi) \subset (M, \varphi)^{\mathcal U}$ is the diagonal inclusion. By \cite[Proposition 3.3]{HI15}, we know that $\rB(M, \varphi) = ((M^{\mathcal U})_{\varphi^{\mathcal U}})' \cap M$. Since $(\rL(\Z), \tau_\Z) \subset ((M^{\mathcal U})_{\varphi^{\mathcal U}}, \varphi^{\mathcal U})$, \cite[Proposition 3.1]{Ue11} implies that 
$$\rB(M, \varphi) \subset \rL(\Z)' \cap M \cap  (M, \varphi) \ast (\rL(\Z), \tau_\Z) = \rL(\Z) \cap M = \C 1.$$

$(\rm iv) \Rightarrow (\rm i)$ When $(M, \varphi)$ is a tracial factor, the implication follows by Popa's result \cite{Po95}. When $\lambda \in (0, 1)$ and $(M, \varphi)$ is a type ${\rm III_\lambda}$ factor endowed with its $\frac{2\pi}{|\log(\lambda)|}$-periodic faithful normal state, the implication follows by \cite[Theorem A]{HI14}. Therefore, we may assume that $M$ is a type ${\rm III_1}$ factor. By \cite[Theorem 3.1]{Ha85}, there exists a faithful normal state $\theta \in M_\ast$ such that $(M_\theta)' \cap M = \C 1$. Then by \cite[Theorem 4]{CS76}, we can choose a sequence of faithful normal states $\varphi_k \in M_\ast$ such that $\lim_k \|\varphi_k - \varphi\| = 0$ and $(M_{\varphi_k})' \cap M = \C 1$ for every $k \in \N$. Choose a nonprincipal ultrafilter $\mathcal U$ on $\N$. By \cite[Theorem A]{HI14}, for every $k \in \N$, there exists a unitary $u_k \in \mathscr U((M^{\mathcal U})_{\varphi_k^{\mathcal U}})$ such that $M$ and $u_k M u_k^*$ are $\ast$-free inside $M^{\mathcal U}$ with respect to $\varphi_k^{\mathcal U}$. Choose another nonprincipal ultrafilter $\mathcal V$ on $\N$ and define $\mathcal W = \mathcal V \otimes \mathcal U$, which is a nonprincipal ultrafilter on $\N \times \N$. We have $(\varphi_k^{\mathcal U})^{\mathcal V} = \varphi^{\mathcal W}$. Set $u = (u_k)^{\mathcal V} \in (M^{\mathcal U})^{\mathcal V}$ and observe that $u \in \mathscr U((M^{\mathcal W})_{\varphi^{\mathcal W}})$. Then $M$ and $u M u^*$ are $\ast$-free inside $M^{\mathcal W}$ with respect to $\varphi^{\mathcal W}$. This further implies that the first factor inclusion $(M, \varphi) \subset (M, \varphi) \ast (M, \varphi)$ is existentially closed and so $(M, \varphi)$ is selfless.
\end{proof}

\begin{remark}
It is well known that all type ${\rm III_1}$ free Araki--Woods factors have trivial bicentralizer \cite{Ho08}. Another proof of this fact was recently obtained in \cite{HI20}. More precisely, in the particular case $q = 0$, the proof of \cite[Main Theorem]{HI20} shows that for any type ${\rm III_1}$ free Araki--Woods factor $\Gamma(H_\R, U)\dpr$ for which the orthogonal representation $U : \R \curvearrowright H_\R$ has a nonzero weakly mixing part, the W$^*$-probability space $(\Gamma(H_\R, U)\dpr, \varphi_U)$ is selfless. It would be interesting to  obtain new classes of type ${\rm III_1}$ factors with trivial bicentralizer by showing that they are selfless.
\end{remark}

\bibliographystyle{plain}

\end{document}